\input amstex
\documentstyle{amsppt}
\topmatter
\NoRunningHeads
\title
On the Weyl asymptotic formula for  Euclidean  domains of finite volume
\endtitle
\author
Leonid Friedlander
\endauthor
\affil
University of Arizona
\endaffil
\email
friedlan@math.arizona.edu
\endemail
\abstract
We give a simple proof of the Weyl asymptotic formula for eigenvalues of the Dirichlet Laplacian, the buckling problem, and the Dirichlet bilaplacian in  Euclidean domain
of finite volume, with no assumption about its boundary.
\endabstract
\endtopmatter
\document

Let $\Omega$ be a domain in $\Bbb R^n$ of finite volume $|\Omega|$. Let $0<\lambda_1\leq \lambda_2\leq\cdots$ be the eigenvalues of the Dirichlet Laplacian in $\Omega$,
and let $0<\mu_1\leq\mu_2\leq\cdots$ be the eigenvalues of the buckling problem in $\Omega$: in the case when the boundary of $\Omega$ is smooth enough, these are 
the values of $\mu$ for which the problem $\Delta^2 v+\mu\Delta v=0$, with the conditions $v=\partial v/\partial\nu=0$ on the boundary of $\Omega$, has a non-trivial solution. Here $\nu$ is the outward 
normal derivative. In general, they are critical values of the functional $\int_\Omega |\Delta v|^2dx$ in $H_0^2(\Omega)$, subject to the constraint
$\int_\Omega |\nabla v|^2dx=1$. We will also consider the Dirichlet bilaplacian. Its eigenvalues, $\omega_j^2$, are critical points of the functional $\int_\Omega |\Delta w|^2dx$ in $H_0^2(\Omega)$, subject to the constraint $\int_\Omega |w(x)|^2dx=1$

Let $N_D(\lambda)$ be the counting function of the sequence $\{\lambda_j\}$: this is the number of $\lambda_j$'s that are strictly less than $\lambda$,  let $N_b(\lambda)$ be the 
counting function of the sequence $\{\mu_j\}$, and let $N_{bl}(\lambda)$ be the counting function of the sequence $\{\omega_j\}$ (note that, for the bilaplacian, it
is the counting function of the square roots of the eigenvalues, not of the eigenvalues themselves.) We will use $N(\lambda)$ when a formula holds for all three counting functions.
It is well known (e.g., see [SV]) that, in the case when $\Omega$ is a bounded domain with smooth boundary,
the Weyl asymptotics holds:
$$\lim_{\lambda\to\infty}\frac{N(\lambda)}{\lambda^{n/2}}=\frac{(4\pi)^{-n/2}}{\Gamma\bigl(\frac{n}{2}+1\bigr)}|\Omega|.\tag 1$$
In this note, I will prove
\proclaim{Theorem}
Let $\Omega$ be an open set in $\Bbb R^n$ of finite volume $|\Omega|$. Then (1) holds.
\endproclaim

This note is motivated by a recent paper [BLPS] where the Weyl asymptotic formula for the buckling problem was proved for domains of finte volume.
The first version of that paper used an additional assumption;  then the authors noticed that it is not needed.  For the Dirichlet Laplacian and the Dirichlet bilaplacian, the results are also  not new. They were proved
by Rozenblum in [R]. He proved the Weyl asymptotic formula for more general class of problems in open sets of finite volume. For the Dirichlet Laplacian, there is another proof
by Simon [S, Theorem 10.6]. Simon makes an assumption that the Lesbegue measure of the bounday is $0$ but, as Michiel van den Berg pointed out to me,
this assumption does not seem to be used in the proof. So, this note contains new proofs of known results; in my view, these proofs are quite simple.

\demo{Proof}
The first ingredient of the proof are well known inequalities
$$N_b(\lambda)\leq N_{bl}(\lambda)\leq N_D(\lambda)\tag 2$$
I give their detailed derivation for the sake of completeness.
They are a consequence of the min-max formulas. Namely
$$N_D(\lambda)=\max\dim \biggl\{\Cal L\subset H_0^1(\Omega):\  \int_\Omega |\nabla u|^2dx<\lambda \int_\Omega |u|^2dx,\ u\in \Cal L \setminus 0\biggr\},$$
$$N_b(\lambda)=\max\dim \biggl\{\Cal L\subset H_0^2(\Omega):\  \int_\Omega |\Delta u|^2dx<\lambda \int_\Omega |\nabla u|^2dx,\ u\in \Cal L \setminus 0\biggr\},$$
$$N_{bl}(\lambda)=\max\dim \biggl\{\Cal L\subset H_0^2(\Omega):\  \int_\Omega |\Delta u|^2dx<\lambda^2 \int_\Omega |u|^2dx,\ u\in \Cal L \setminus 0\biggr\}.$$
 To show that $ N_{bl}(\lambda)\leq N_D(\lambda)$ we notice that if
$u(x)\in H_0^2(\Omega)\subset  H^2(\Omega)\cap H_0^1(\Omega)$ and
$$ \int_\Omega |\Delta u|^2dx<\lambda^2 \int_\Omega | u|^2dx$$
then
$$\int_\Omega |\nabla u|^2dx=-\int_\Omega \Delta u\cdot \bar u dx\leq \lambda
\int_\Omega | u|^2dx;$$
here we used the Cauchy-Schwarz inequality. To verify the inequality $N_b(\lambda)\leq N_{bl}(\lambda)$, let us take
$u(x)\in H_0^2(\lambda)$ such that 
$$\int_\Omega |\Delta u|^2dx<\lambda \int_\Omega |\nabla u|^2dx.$$
Then
$$\int_\Omega |\Delta u|^2dx<\lambda \int_\Omega |\nabla u|^2dx=\lambda\int_\Omega \Delta u\cdot \bar u dx\leq\lambda
\biggl(\int_\Omega |\Delta u|^2dx\biggr)^{1/2}\biggl(\int_\Omega | u|^2dx\biggr)^{1/2};$$
this implies
$$ \int_\Omega |\Delta u|^2dx<\lambda^2 \int_\Omega |u|^2dx.$$

The second well known fact that follows immediately from the min-max formulas is that  all three counting functions are superadditive with respect to the domain: if
$\Omega_1$, $\Omega_2$ are open sets,  $\Omega_1\subset\Omega$, $\Omega_2\subset\Omega$, and $\Omega_1\cap\Omega_2=\emptyset$
then ([CH], VI.2, Theorem 2)
$$N(\lambda,\Omega)\geq N(\lambda, \Omega_1)+N(\lambda,\Omega_2). \tag 3$$

Property (3) implies
$$\liminf_{\lambda\to\infty}\frac{N(\lambda)}{\lambda^{n/2}}\geq\frac{(4\pi)^{-n/2}}{\Gamma\bigl(\frac{n}{2}+1\bigr)}|\Omega|.\tag 4$$
To prove (4) we will show that for every $\epsilon>0$ there exist a finite number of bounded domains $\Omega_j\subset \Omega$ with smooth boundary that are mutually disjoint,
and that their total volume is bigger that $|\Omega|-\epsilon$. Then
$$\liminf_{\lambda\to\infty}\frac{N(\lambda,\Omega)}{\lambda^{n/2}}\geq
\sum_j \lim_{\lambda\to\infty}\frac{N(\lambda,\Omega_j)}{\lambda^{n/2}}\geq
\frac{(4\pi)^{-n/2}}{\Gamma\bigl(\frac{n}{2}+1\bigr)}(|\Omega|-\epsilon),$$
which implies (4). I used the fact that the Weyl asymptotics holds for bounded domains with smooth boundary.
To construct $\Omega_j$, we introduce $\Omega_\eta=\{x\in\Omega: \text{dist}(x,\Bbb R^n\setminus\Omega)>\eta\}$; here $\eta$ is a positive number  By the dominated convergence theorem,
$\lim_{\eta\to 0}|\Omega_\eta|=|\Omega|$. Find $\eta$ such that $|\Omega_\eta|>|\Omega|-(\epsilon/2)$, cover $\Bbb R^n$ by cubes of side $\eta/\sqrt{n}$, pick up the cubes 
that lie completely in $\Omega$, and, inside of each one take a smooth subdomain of the volume greater than the volume of a cube minus $\epsilon$ divided by two times the number of cubes.

To prove the Weyl asymptotics for $N_D(\lambda)$, I will use the heat trace
$$h(t)=\int_0^\infty e^{-t\lambda}dN_D(\lambda)=t\int_0^\infty e^{-t\lambda}N_D(\lambda)d\lambda.$$
Property (4) easily implies
$$\liminf_{t\to 0}t^{n/2}h(t)\geq (4\pi)^{-n/2}|\Omega|.\tag 5$$
Indeed, for every $C$,
$$C<\frac{(4\pi)^{-n/2}}{\Gamma\bigl(\frac{n}{2}+1\bigr)}|\Omega|,$$
there exists $a>0$ such that $N_D(\lambda)\geq C\lambda^{n/2}$ when $\lambda\geq a$.
Then
$$h(t)\geq Ct\int_a^\infty \lambda^{n/2}e^{-t\lambda}d\lambda =Ct^{-n/2}\int_{at}^\infty \lambda^{n/2}e^{-\lambda}d\lambda.$$

Let $H_D(x,y,t)$ be the heat kernel of the Dirichlet Laplacian in $\Omega$. By the domain monotonicity of Dirichlet heat kernel,
$$H_D(x,y,t)\leq (4\pi t)^{-n/2}\exp\biggl(\frac{|x-y|^2}{4t}\biggr);$$
here the right hand side  is the heat kernel in $\Bbb R^n$. Therefore,
$$h(t)=\int_\Omega H_D(x,x,t)dx\leq (4\pi t)^{-n/2}|\Omega|.$$
The last inequality and  (5) give us
$$\lim_{t\to 0}t^{n/2}h(t)= (4\pi)^{-n/2}|\Omega|.$$
Then, by the Karamata tauberian theorem, (1) holds for $N_D(\lambda)$.

The Weyl asymptotics for $N_D(\lambda)$, together with (2), implies
$$\limsup_{\lambda\to\infty}\frac{N_b(\lambda)}{\lambda^{n/2}}\leq\frac{(4\pi)^{-n/2}}{\Gamma\bigl(\frac{n}{2}+1\bigr)}|\Omega|, \quad
\limsup_{\lambda\to\infty}\frac{N_{bl}(\lambda)}{\lambda^{n/2}}\leq\frac{(4\pi)^{-n/2}}{\Gamma\bigl(\frac{n}{2}+1\bigr)}|\Omega|.\tag 6$$
The statement of the theorem follows from  (4) and (6).

\enddemo
\proclaim{Acknowledgements}
I would like to thank Michael Levitin and Michiel van den Berg for reading the manuscript and giving important feedback.
Michiel van den Berg pointed my attention to [R] and [S].
\endproclaim

\enddocument